\documentstyle[12pt,twoside]{amsart}

\newtheorem{thm}{Theorem}[section]

\newtheorem{lem}[thm]{Lemma}
\newtheorem{cor}[thm]{Corollary}

\newtheorem{prop}[thm]{Proposition}

\newtheorem{complement}[thm]{Complement}

\theoremstyle{definition}
\newtheorem{defn}[thm]{Definition}

\newtheorem{say}[thm]{}
\newtheorem{exmp}[thm]{Example}

\newtheorem{const}[thm]{Construction}   

\newtheorem{rem}[thm]{Remark}          
\newtheorem{ack}{Acknowledgments}        
\newtheorem{notation}[thm]{Notation}

\theoremstyle{remark}

\renewcommand{\o}[0]{{\cal O}}
\renewcommand{\c}[0]{{\Bbb C}}   
\newcommand{\z}[0]{{\Bbb Z}}
\newcommand{\n}[0]{{\Bbb N}}
\newcommand{\r}[0]{{\Bbb R}}
\newcommand{\p}[0]{{\Bbb P}}

\newcommand{\map}[0]{\dasharrow}
\newcommand{\qtq}[1]{\quad\mbox{#1}\quad}
\newcommand{\spec}[0]{\operatorname{Spec}}

\newcommand{\mult}[0]{\operatorname{mult}}

\newcommand{\Hom}[0]{\operatorname{Hom}}

\newsymbol\subsetneq 2328
\newsymbol\onto 1310
\newsymbol\into 1323
\newsymbol\twoheadrightarrow 1310

\begin{document}
\bibliographystyle{amsplain}

\title{Real Algebraic Threefolds I.\\ Terminal Singularities}
\author{J\'anos Koll\'ar}

\maketitle
\tableofcontents

\section{Introduction}

In real algebraic   geometry, considerable attention has
been paid to the study of real algebraic curves (in
connection with Hilbert's 16th problem) and also to real
algebraic surfaces. See \cite{Viro90}, \cite{Riesler93} and the
references there.

In higher dimensions one of the main avenues of
investigation was initiated by \cite{Nash52},  and later
developed by many others (see
\cite{AK92} for some recent directions). One of these results says
that every compact differentiable manifold can be realized
as the set of real points of an algebraic variety.
\cite{Nash52} posed the problem of obtaining similar results
using a restricted class of varieties, for instance rational
varieties.

The aim of this series of papers is to develop the  theory of 
minimal models for
real algebraic threefolds. This approach gives very strong
information about the topology of real algebraic threefolds, and 
 it  also answers the above mentioned question of
\cite{Nash52}.

For algebraic threefolds over $\c$, the minimal model
program provides a very powerful tool.  The method of the
program is the following. (See \cite{Koll87} or \cite{CKM88}
for introductions)

Starting with a smooth projective variety $X$, we perform a
series of ``elementary" birational transformations
$$ X=X_0\map X_1\map \cdots \map X_n
$$ until we reach a variety $X_n$ whose global structure is
``simple". In essence the minimal model
program  allows us to investigate many questions in two
steps: first study the effect of the ``elementary"
transformations and then consider the ``simple" global
situation.

In parctice both of these steps are frequently rather
difficult. For instance, we still do not have a complete
list of all possible ``elementary" steps, despite repeated
attempts to obtain it.

 A somewhat unpleasant feature of the theory is
that the varieties $X_i$ are not smooth, but have so called
terminal singularities. In developing the theory of minimal models
for real algebraic threefolds, we again have to understand the
occurring terminal singularities.

The aim of this paper is to give  a classification of 
terminal 3-fold singularities over
$\r$. Minimal models serve only as a background, the proofs depend
entirely on well established methods of singularity theory. 
I do not even  use the definition of terminal singularities!

Terminal 3-fold singularities over $\c$ are completely
classified.
\cite{Reid85} is a very readable   introduction and survey.
I will take the result of this classification as my definition,
since the theory over $\r$ can be most naturally developed in this
setting.

The classification is, in some sense, not complete.
 In a few cases I obtain unique normal forms
(\ref{ca1-top.thm}), but in most cases this seems nearly
impossible (see
\cite{Markushevich85} for a
 special case over $\c$). My aim is to write the
singularities in a form that allows one to determine their
topology over $\r$. The resulting lists and algorithms are
given in sections 4--5.

It turns out that the normal forms of 3-fold terminal
singularities are essentially the same over any field of
characteristic zero. Thus in sections 2--3 I work with any
subfield of $\c$. 

As a consequence of the classification over $\c$, we know
that  3-fold terminal singularities come in two types. Some
are hypersurface singularities, and the others are quotients
of these hypersurface singularities by a finite cyclic
group. Accordingly, the classification over any field is
 done in two steps. Section 2   deals with  terminal
hypersurface singularities. These results are mostly routine
generalizations of the theory over $\c$.

Quotient singularities frequently have ``twisted" forms over
a subfield of $\c$.   ``Twisted" forms do not
appear for  3-fold terminal singularities, and so the
classification ends up very similar to the one over $\c$. 

\begin{ack}  I   thank G. Mikhalkin for
answering my numerous questions about  real algebraic geometry. 
Partial financial support was provided by  the NSF under grant
number  DMS-9622394. Most of this paper was written while I visited
RIMS, Kyoto Univ.
\end{ack}

\section{Terminal hypersurface singularities}

\begin{notation} For a field $K$ let $K[[x_1,\dots,x_n]]$
denote the ring of formal power series in $n$ variables over
$K$. For $K=\r$ or
$K=\c$, let $K\{x_1,\dots,x_n\}$ denote the ring of those
formal power series which converge in some  neighborhood of
the origin.

For any $F\in K\{x_1,\dots,x_n\}$  the set $(F=0)$ is a germ
of a real or complex analytic set. I will refer to it as a
singularity. If $F\in K[[x_1,\dots,x_n]]$ then by the
singularity $(F=0)$ I mean the scheme
$\spec_KK[[x_1,\dots,x_n]]/(F)$. 

For a power series $F$, $F_d$ denotes the degree $d$
homogeneous part. The multiplicity, denoted by $\mult_0F$,
is the smallest $d$ such that
$F_d\neq 0$. If we write a power series as $F_{\geq d}$ then
it is assumed that its multiplicity is at least $d$.

Two power series $F,G\in K[[x_1,\dots,x_n]]$
are called equivalent over $K$ if there is an automorphism  of
$K[[x_1,\dots,x_n]]$ given by
$x_i\mapsto \phi_i(x_1,\dots,x_n)\in  K[[x_1,\dots,x_n]]$
and an invertible $u(x_1,\dots,x_n)\in K[[x_1,\dots,x_n]]$
such that
$$
u(x_1,\dots,x_n)G(x_1,\dots,x_n)=F(\phi_1,\dots,\phi_n).
$$ Thus $F$ and $G$ are equivalent iff the corresponding
singularities $(F=0)$ and $(G=0)$ are isomorphic (over
$K$).

We  have to pay special attention to cases when $F$
and $G$ are not equivalent over $K$ but are equivalent
over some larger field. For instance, $F=x_1^2+x_2^2$ and
$G=x_1^2-x_2^2$ are not equivalent over  $\r$ but are
equivalent over $\c$.

If $K=\r,\c$ and $F,G\in K\{x_1,\dots,x_n\}$ then I am
mainly interested in equivalences where $u,\phi_i\in
K\{x_1,\dots,x_n\}$.

If $F,G\in K\{x_1,\dots,x_n\}$ have  isolated critical
points at the origin, then
$F$ and $G$ are equivalent in $K\{x_1,\dots,x_n\}$ iff
they are equivalent in $K[[x_1,\dots,x_n]]$ (cf.
\cite[p.121]{AGV85}), thus we do not have to be  careful
about this distinction.
\end{notation}

\begin{defn} Let $K$ be a field of characteristic zero with
algebraic closure $\bar K$. $(F(x,y,z)=0)$ is called a {\it Du
Val} singularity (or a rational double point) iff over $\bar
K$ it is  equivalent to  one of the  standard forms
\begin{enumerate}
\item[$A_n$] \quad $x^2+y^2+z^{n+1}=0$;
\item[$D_n$] \quad  $x^2+y^2z+z^{n-1}=0$;
\item[$E_6$]  \quad $x^2+y^3+z^4=0$;
\item[$E_7$]  \quad $x^2+y^3+yz^3=0$;
\item[$E_8$]  \quad $x^2+y^3+z^5=0$.
\end{enumerate} Du Val singularities have many interesting
intrinsic characterizations, (cf. \cite{Durfee79, Reid85})
but I will not use this.
\end{defn}

The following  definition introduces our basic objects of study.

\begin{defn}\label{cdv.def} Let $K$ be a field of
characteristic zero with algebraic closure $\bar K$.
$(F(x,y,z,t)=0)$ is called a {\it compound Du Val} singularity (or
{\it cDV} for short) iff over $\bar K$ it is  equivalent to  
$$ h(x,y,z)+tf(x,y,z,t)=0
$$ where $(h=0)$ is a  Du Val singularity.

$(F(x,y,z,t)=0)$ is called a
$cA_n$ (resp. $cD_n$ or $cE_n$) singularity if its equation
can be written as above with $h$ having type 
$A_n$ (resp. $D_n$ or $E_n$), but it does not admit such
representation with a smaller value of $n$.  It is called a
$cA$ (resp. $cD$ or $cE$) singularity if the value of $n$ is
not specified.
\end{defn}

The reason  we are interested in cDV singularities is the
following:

\begin{thm}\cite{Reid80} A 3-dimensional hypersurface
singularity over
$\c$ is terminal iff it is an isolated cDV singularity.\qed
\end{thm}

The aim of this section is to develop ``normal forms" for cDV
singularities over any field $K$. This will then give 
``normal forms" for 3-dimensional terminal hypersurface
singularities over $K$.

The proof is a rather standard application of the
methods  of
\cite{AGV85}. 

\begin{say}\label{nf.meth} We repeatedly use 3 methods:
\begin{enumerate}
\item The Weierstrass preparation theorem. This is
frequently stated only over
$\c$, but it works over any field since the Weierstrass
normal form is unique.
\item The elimination of the $y^{n-1}$-term from the
polynomial
$a_ny^n+a_{n-1}y^{n-1}+\dots$ by a coordinate change
$y\mapsto y-a_{n-1}/na_n$ when $a_n$ is invertible.
\item Let $M_1,\dots,M_k$ be 
monomials in the variables
$x_1,\dots,x_m$. Assume that 
$x_0M_1,\dots,x_0M_k$ are  multiplicatively independent.
Then any power series of the form
$\sum M_i\cdot u_i(x_1,\dots,x_m)$  where $u_i(0)\neq 0$ 
for all $i$  is equivalent to
$\sum M_i\cdot u_i(0)$  by a suitable coordinate change 
$x_i\mapsto x_i\cdot(\mbox{unit})$. 
\end{enumerate}
\end{say}

These elementary operations are sufficient  to deal with the
$cA$ and
$cE$ cases. In the $cD$ case the following  generalization
of  (\ref{nf.meth}.2)  is needed.

\begin{const}\label{nf.meth.weighted}  In
$K[[x_1,\dots,x_m]]$, assign positive integral weights to
the variables $w(x_i)=w_i$. For a monomial set $w(\prod
x_i^{c_i})=\sum c_iw_i$.  Write a power series in terms of
its weighted homogeneous pieces
$F=F_d+F_{d+1}+\dots$. Choose $g_i\in K[[x_1,\dots,x_m]]$
such that
$w(g_i)=w(x_i)+e$ for some $e>0$. Then
$$ F(x_i+g_i)=F(x_i)+\sum g_i\frac{\partial F_d}{\partial
x_i}+ R_{>(d+e)}(x_i).
$$ Repeatedly using this for higher and higher  degrees, we
see that,  for every $N>0$,
$F$ is equivalent to a power series
$F^N+R_{>N}$ where $F^N$ is a polynomial of degree $N$ and
no linear combination of the monomials in $F^N$ can be
written in the form
$\sum g_i(\partial F_d/\partial x_i)$ as above.

In the ring of formal power series this can be continued
indefinitely, thus at the end we can kill all the degree
$>d$ elements of the  Jacobian ideal
$$
\Delta(F_d):=\left(\frac{\partial F_d}{\partial x_1},
\dots , \frac{\partial F_d}{\partial x_m}\right).
$$

If $F\in K\{x_1,\dots,x_m\}$  defines an isolated
singularity, then by Tougeron's lemma (cf.
\cite[p.121]{AGV85}), $F^N+R_{>N}$ is equivalent to
$F^N$ by an analytic coordinate change for $N\gg 1$. Thus
the final conclusion is the same.
\end{const}

\begin{prop}\label{morse.lem} Any power series $F_{\geq
2}(x_1,\dots,x_n)$ is equivalent to a power series
$$ a_1x_1^2+\dots+a_kx_k^2+G_{\geq 3}(x_{k+1},\dots,x_n).
$$
\end{prop}

Proof. By a linear change of coordinates  we can diagonalize
$F_2$, thus we can assume that
$F_2=a_1x_1^2+\dots+a_kx_k^2$. Repeatedly applying
(\ref{nf.meth}.1) to the variables
$x_1,\dots,x_k$ we reach a situation when $F$ is a quadratic
polynomial in
 the variables $x_1,\dots,x_k$. (\ref{nf.meth}.2) can then
be used to eliminate the linear terms in $x_1,\dots,x_k$.\qed

\begin{thm}\label{ca.thm} Assume that  $F_{\geq
1}(x,y,z,t)\in K[[x,y,z,t]]$ defines a terminal singularity
of type $cA$. Then $F$ is equivalent to one of the following:
\begin{enumerate}
\item[$cA_0$] \quad $x=0$.
\item[$cA_1$]  \quad $ax^2+by^2+cz^2+dt^m=0$, where
$abcd\neq 0$.
\item[$cA_{>1}$]  \quad $ax^2+by^2+f_{\geq 3}(z,t)=0$, where
$ab\neq 0$. This has type $cA_n$ for $n=\mult_0f-1$. 
\end{enumerate}
\end{thm}

Proof. If $F_1\neq 0$ then (\ref{nf.meth}.1) gives $cA_0$.
Thus assume that $F_1=0$. $F$ has type $cA$, hence $F_2$ is
a quadric of rank at least 2. If the rank is 2 then
(\ref{morse.lem}) gives the 
$cA_{>1}$ cases.

Assume finally that $F_2$ has rank 3 or 4.
 By (\ref{morse.lem}) we can write $F$ as
$ax^2+by^2+cz^2+g(t)=0$. Using (\ref{nf.meth}.3) we obtain
$ax^2+by^2+cz^2+dt^m=0$.

In all these cases we can multiply through by $a^{-1}$ to get a
somewhat simpler form when the coefficient of $x^2$ is 1.
\qed

\begin{thm}\label{cd.thm} Assume that  $F_{\geq
2}(x,y,z,t)\in K[[x,y,z,t]]$ defines a terminal singularity
of type $cD$. Then $F$ is equivalent to one of the following:
\begin{enumerate}
\item[$cD_4$] \quad $x^2+f_{\geq 3}(y,z,t)$, where $f_3$ is
not divisible by the square of a linear form.
\item[$cD_{>4}$] \quad $x^2+y^2z+ayt^r+h_{\geq s}(z,t)$,
where $a\in K$,
 $r\geq 3$, $s\geq 4$ and $h_s\neq 0$. This has type  $cD_n$
where
$n=\min\{2r, s+1\}$ if $a\neq 0$ and $n=s+1$ if $a=0$. 
\end{enumerate}
\end{thm}

Proof. $F_2$ is a rank one quadric, thus in suitable
coordinates the equation becomes  $ax^2+f_{\geq 3}(y,z,t)$.
Here $f_3\neq 0$  is not the cube of a linear form since
otherwise we would have a type $cE$ singularity.  If $f_3$
is not divisible by the square of a linear form then we have
case
$cD_4$.

If $f_3$ is  divisible by the square of a linear form, then
$f_3=l_1^2l_2$ for two linear forms $l_i$, and both of them
are defined over $K$. We can change coordinates  $l_1\mapsto
y$ and $l_2\mapsto z$.

At this point our power series is $x^2+y^2z+(\mbox{higher
order terms})$. Assign weights $w(x)=3, w(y)=w(z)=2,
w(t)=6$. The leading term is
$x^2+y^2z$. Using (\ref{nf.meth.weighted}) we can eliminate
all monomials which contain $y^2$ or $yz$.

To see the last part, take the hyperplane section
$t=\lambda z$. The term
$ay\lambda^rz^r$ can be eliminated by a substitution
$y\mapsto y+(a/2)\lambda^rz^{r-1}$. This creates  a term
$-(a/2)^2\lambda^{2r} z^{2r-1}$. The only problem could be
that
$h(z,\lambda z)$ has multiplicity $2r-1$ and there is
cancellation. However, $h_{2r-1}(z,\lambda z)=z^{2r-1}
h_{2r-1}(1,\lambda)$ is a polynomial of degree $2r-1$ in
$\lambda$, thus it does not equal
$-(a/2)^2\lambda^{2r} z^{2r-1}$.
\qed

\begin{thm}\label{ce.thm} Assume that  $F_{\geq
2}(x,y,z,t)\in K[[x,y,z,t]]$ defines a terminal singularity
of type $cE$. Then $F$ is equivalent to one of the following:
\begin{enumerate}
\item[$cE_6$] \quad $x^2+y^3+yg_{\geq 3}(z,t)+h_{\geq
4}(z,t)$, where
$h_4\neq 0$.
\item[$cE_7$] \quad $x^2+y^3+yg_{\geq 3}(z,t)+h_{\geq
5}(z,t)$, where
$g_3\neq 0$.
\item[$cE_8$] \quad $x^2+y^3+yg_{\geq 4}(z,t)+h_{\geq
5}(z,t)$, where
$h_5\neq 0$.
\end{enumerate}
\end{thm}

Proof. $F_2$ is a rank one quadric by (\ref{cdv.def}), thus
in suitable coordinates the equation becomes  $ax^2+f_{\geq
3}(y,z,t)$. Here
$f_3\neq 0$ and it is the cube of a linear form since
otherwise we would have a type $cD$ singularity. 
(\ref{nf.meth}.1--2) gives an equation
$$ ax^2+by^3+yg_{\geq 3}(z,t)+h_{\geq 4}(z,t).
$$ Multiply the equation by $a^3b^2$ and then make the
substitutions
$x\mapsto xa^{-2}b^{-1}$ and $y\mapsto ya^{-1}b^{-1}$ to
get the required normal forms.\qed

\section{Higher index terminal singularities}

The classification of non-hypersurface terminal 3-fold
singularities over $\c$ relies on the following construction:

Let $\z_n$ denote the cyclic group of order $n$ and
$\epsilon$ a primitive
$n^{th}$ root of unity. Assume that $\z_n$ acts on $\c^4$ by
$$
\sigma: (x,y,z,t)\mapsto 
(\epsilon^{a_x}x,\epsilon^{a_y}y,\epsilon^{a_z}z,\epsilon^{a_t}t).
$$ I will use the shorter notation
$\frac{1}{n}(a_x,a_y,a_z,a_t)$ to denote such an action.

If $F(x,y,z,t)$ is equivariant with respect to this action,
then
$\z_n$ acts on the hypersurface $(F=0)$ and we can take the
quotient,
 denoted by $(F=0)/\frac{1}{n}(a_x,a_y,a_z,a_t)$.

By \cite{Reid80},  every terminal 3-fold
singularity $X$ over
$\c$ is of the form
$(F=0)/\frac{1}{n}(a_x,a_y,a_z,a_t)$, where $F$ defines a
terminal hypersurface singularity. The value of $n$ is
uniquely determined by
$X$, it is called the {\it index} of $X$.

It is not easy to come up with a complete list of  terminal
3-fold singularities, but by now the list is well
understood; see
\cite{Reid85} for  a good survey.  It turns out that most
actions do not produce terminal quotients and we have only a
few cases:

\begin{thm}\label{hind.overC.thm}\cite{Mori85}
 Let $0\in X$ be a 3-fold terminal
nonhypersurface singularity over
$\c$. Then $0\in X$ is isomorphic   to a singularity
described by the following list:
$$
\begin{tabular}{|c|l|c|l|l|}
\hline name &\qquad equation  & index  & \quad action &
condition\\
\hline cA/n& $xy+f(z,t)$ & $n$ & $(r,-r,1,0)$ & $(n,r)=1$\\
\hline cAx/2&  $x^2+ y^2 +f_{\geq 4}(z,t)$  & $2$ & $(0,1,1,1)$&\\
\hline cAx/4 & $x^2+y^2 +f_{\geq 2}(z,t)$ & $4$ & $(1,3,1,2)$&
$f_2(0,1)=0$\\
\hline cD/2&$x^2+f_{\geq 3}(y,z,t)$ & $2$ &$(1,0,1,1)$&\\
\hline cD/3&$x^2+f_{\geq 3}(y,z,t)$ &$3$ &$(0,2,1,1)$&
            $f_3(1,0,0)\neq 0$\\
\hline cE/2&$x^2+y^3+f_{\geq 4}(y,z,t)$ &$2$
&$(1,0,1,1)$&\\
\hline
\end{tabular}
$$ 
\end{thm}

The equations have to satisfy 2 obvious conditions: 
\begin{enumerate}
\item The equations   define a terminal hypersurface
singularity.

\item The equations   are $\z_n$-equivariant. (In fact
$\z_n$-invariant, except for $cAx/4$.)
\end{enumerate}

If we work over a field $K$ which does not contain
the
$n^{th}$ roots of unity, then the action
$\frac{1}{n}(a_1,\dots,a_m)$ is not defined over $K$. There
is, however, another way of loking at the quotient which
does make sense over any field.

Any action of the cyclic group $\z_n$ on $\c^m$ defines a
$\z_n$-grading
$w$ of $\c[[x_1,\dots,x_m]]$  by
$$ w(\prod x_i^{c_i})=a \qtq{iff} 
\sigma (\prod x_i^{c_i})=\epsilon^a \cdot \prod x_i^{c_i}.
$$ If $F$ is $\z_n$-equivariant then 
$(F)\subset \c[[x_1,\dots,x_m]]$ is a homogeneous ideal,
hence the grading descends to a grading of
$\c[[x_1,\dots,x_m]]/(F)$. The ring of  functions  on the
quotient
$(F=0)/\frac{1}{n}(a_1,\dots,a_m)$ can be identified with
the ring of grade zero elements of $\c[[x_1,\dots,x_m]]/(F)$.

If $K$ is any field, $n\in \n$  and  $a_i\in \z$, then we
obtain a
$\z_n$-grading $w=w(a_1,\dots,a_m)$ of $K[[x_1,\dots,x_m]]$ 
(or of
$\r\{x_1,\dots,x_m\}$) by
$$ w(\prod x_i^{c_i})=\sum c_ia_i \in \z_n.
$$ Let $R\subset K[[x_1,\dots,x_m]]$  denote the subring of
grade zero elements.  Then $\spec_KR$ gives a singularity
over $K$ which is denoted by
$$ {\Bbb A}^m/{\textstyle \frac{1}{n}}(a_1,\dots,a_m).
$$ (Especially when $K=\r$, one might be tempted to write
${\Bbb R}^m/{\textstyle \frac{1}{n}}(a_1,\dots,a_m)$
instead. However,  
the set of real points of ${\Bbb A}^m/{\textstyle
\frac{1}{n}}(a_1,\dots,a_m)$   is not in any sense a
quotient of the  set $\r^n$ (cf. (\ref{2ind.quot})),  so this may
lead to confusion.)

If $F\in K[[x_1,\dots,x_m]]$ is graded homogeneous, then
$w$ gives a grading of $K[[x_1,\dots,x_m]]/(F)$. Let
$R/(R\cap(F))\subset K[[x_1,\dots,x_m]]/(F)$  be the subring
of grade zero elements. 
$\spec_K R/(R\cap(F))$  defines a singularity over $K$. By
construction, 
$$
\spec_K R/(R\cap(F))\times_{\spec K}\spec \bar K\cong 
(F=0)/{\textstyle
\frac{1}{n}}(a_1,\dots,a_m).
$$ Thus $\spec_K R$ is a terminal singularity over $K$  iff 
$(F=0)/{\textstyle \frac{1}{n}}(a_1,\dots,a_m)$ is a 
terminal singularity over $\bar K$.

Under certain conditions,  every $K$-form of a quotient
is obtained this way:

\begin{thm}\label{hind.gen.thm}
 $K$ be a field of characteristic zero with algebraic
closure $\bar K$.
Let $\z_n$ denote the cyclic group of order $n$ and
$\epsilon$ a primitive
$n^{th}$ root of unity. Assume that $\z_n$ acts on $\bar K^m$
by
$\sigma: (x_i)\mapsto 
(\epsilon^{w_i}x_i)$.
 Let $F\in \bar K[[x_1,\dots,x_m]]$ be equivariant with respect to this
action, and assume that the fixed point set of $\sigma$ has
codimension at least 2 in $(F=0)$. 
Assume in addition that
$$
w(F)-\sum w_i\qtq{is relatively prime to} n.
$$
Let  $0\in X$ be a singularity over $K$
such that 
$$
X\times_{\spec K}\spec \bar K\cong (F=0)/
\textstyle{\frac{1}{n}(w_1,\dots,w_m)}.
$$
Then there is an $F^K\in K[[x_1,\dots,x_m]]$ such that
$F$ and $F^K$ are equivalent over $\bar K$ and
$$
X\cong (F^K=0)/
\textstyle{\frac{1}{n}(w_1,\dots,w_m)}.
$$
\end{thm}

It is worthwhile to note that the condition about
$n$ and $w(F)-\sum w_i $ being  relatively prime  is essential:

\begin{exmp} Consider the  quotient singularity
 $\c[u,v]/\frac{1}{n}(1,-1)$. It is isomorphic to
$(xy-z^n=0)$ via the
substitutions $x=u^n,y=v^n,z=uv$. Over $\c$  we have a 
 Du Val singularity
$A_{n-1}=(x^2+y^2+z^n=0)$. 

Over $\r$ we see that  
$(x^2-y^2-z^n=0)\cong {\Bbb A}^2/\frac{1}{n}(1,-1)$. Another
$\r$-form of $A_{n-1}$ is $x^2+y^2-z^n$. This can also be
obtained as a quotient, but this time we act on ${\Bbb A}^2$
by rotation with angle
$2\pi/n$.

Finally, if $n$ is even, then there is another $\r$-form of
$A_{n-1}$ given by
$(x^2+y^2+z^n=0)$.  The only $\r$-point is the origin, so we
do not even have a nonzero map
$\r^2\to (x^2+y^2+z^n=0)$. 

As  another example, take the 4-dimensional terminal
singularity
$\c^4/\frac{1}{n}(a,-a,b,-b)$ for any $(ab,n)=1$.  It has
another
$\r$-form given as ${\Bbb A}^4/\z_n$ where we act on the
first two coordinates by rotation with angle $2a\pi/n$ and
on  the last two coordinates by rotation with angle
$2b\pi/n$.

In some special cases there are further $\r$-forms. Take for
instance
$\c^4/\frac{1}{2}(1,1,1,1)$. This can be realized as the
cone over
$\c\p^3$ embedded by the quadrics to $\c\p^9$. 

Let $C\subset \r\p^2$ be a smooth conic. Taking symmetric powers
we have $S^3C\subset S^3\r\p^2$ and $S^3H^0(\r\p^2,\o(1))$ embeds
it to $\r\p^9$. If $C$ has a real point, then 
$S^3C\cong \r\p^3$ and we get the Veronese embedding.  If $C$ has no
real points then the image is a variety   over $\r$ without real
points. The cone over it  is  a real form of
${\Bbb A}^4/\frac12(1,1,1,1)$ with an  isolated real point at the
origin.
\end{exmp}

Proof of (\ref{hind.gen.thm}). 
Set $S=\bar K[[x_1,\dots,x_m]]/(F)$,
the ring of functions on $\tilde X_{\bar K}:=(F=0)$. 
 The $\z_n$-action defines a
$\z_n$-grading   $S=\sum_{i=0}^{n-1}S_i$. $S_0$, the ring of grade $0$
elements,  is exactly  the ring of functions on $X_{\bar K}$.  Our aim is
to find an algebraic way of reconstructing $S$ from $S_0$, which then
hopefully  generalizes to nonclosed fields.

There is another summand  which can be easily seen algebraically.
Set $d=w(F)-\sum w_i$.  Note that
$$
\frac{1}{\partial F/\partial x_m}dx_1\wedge\dots\wedge dx_{m-1}
$$
is a local generator of $\omega_S$ and it has weight $-d$. Thus
$$
\omega_{S_0}\cong S_d\frac{1}{\partial F/\partial
x_m}dx_1\wedge\dots\wedge dx_{m-1}.
$$
Once $S_d$ is determined, we obtain $S_{jd}$ 
as follows. The multiplication map
$$
S_a\otimes_{S_0} S_b\to S_{a+b}\qtq{(subscripts modulo $n$)}
$$ are
isomorhisms over the open set where the $\z_n$-action is free. 
We assumed that the complement has codimension at least 2, thus
$S_{jd}\cong S_d^{[j]}$, where  $ S_d^{[j]}$ denotes the double dual
of
$S_d^{\otimes j}$. If $d$ and $n$ are
 relatively prime, then we obtain  every summand $S_i$ this way. In
particular, 
$$
S=\sum_{i=0}^{n-1}S_i\cong \sum_{j=0}^{n-1}\omega_{S_0}^{[j]}.
$$

Over an arbitrary field, we can thus proceed as follows.
Let $\omega_X$ be the dualizing sheaf of $X$.  This is also
the reflexive sheaf $\o_X(K_X)$ where $K_X$ is the canonical
class.

Then $\omega_X^{[n]}$ is isomorphic to $\o_X$, where $n$ is
the index. (We know this over $\bar K$. Isomorphism of two
sheaves
$F,G$ is a question about $\Hom(F,G)$ and this commutes with
base field extensions.) Fix such an isomorphism
$s:\omega_X^{[n]}\to \o_X$.

Consider the $\o_X$-algebra
$$
 R(X,s):=\sum_{j=0}^{n-1}\omega_X^{[j]},
$$ where multiplication for $j+k\geq n$ is given by
$$
\omega_X^{[j]}\otimes \omega_X^{[k]}\mapsto
\omega_X^{[j+k]}\cong \omega_X^{[n]}\otimes
\omega_X^{[j+k-n]}\stackrel{s\otimes 1}{\longrightarrow} 
\omega_X^{[j+k-n]}.
$$ This has a $\z_n$ grading by declaring $\omega_X^{[j]}$
to have grade
$j$. 

(Note. Two   isomorphisms $s_1,s_2:\omega_X^{[r]}\to \o_X$
differ by an invertible function $h\in \o_X^*$. If $h$ is an
$n^{th}$-power, then the resulting algebras $R(X,s_i)$
are isomorphic, but they need not be isomorphic otherwise.
This is connected with the topological aspects observed in
(\ref{2ind.quot}).)

Over $\bar K$, $R(X,s)$ is isomorphic to  $\o_{\tilde X}$. Thus 
$R(X,s)$ is a $K$-form of
$\o_{\tilde X}$. In particular, $R(X,s)$ is an algebra
of the form
$K[[x_1,\dots,x_m]]/(F^K)$, where $F$ and $F^K$ are equivalent over
$\bar K$.

 The grading lifts to a grading of
$K[[x_1,\dots,x_m]]$ such that $F^K$ is graded homogeneous. We can
choose
$x_i$ to be homogeneous.\qed

As a corollary,  we obtain the following classification of
terminal 3-fold nonhypersurface singularities over nonclosed
fields:

\begin{thm}\label{hind.term.thm} Let $K$ be  a field of
characteristic zero and $0\in X$ a 3-fold terminal
nonhypersurface singularity over
$K$. Then $0\in X$ is isomorphic over $K$ to a singularity
described by the following list:
$$
\begin{tabular}{|c|l|c|l|l|}
\hline name &\qquad equation  & index  & \quad weights &
condition\\
\hline cA/2& $ax^2+by^2+f(z,t)$ & $2$ & $(1,1,1,0)$&\\
\hline cA/n& $xy+f(z,t)$ & $n\geq 3$ & $(r,-r,1,0)$&
$(n,r)=1$\\
\hline cAx/2&  $ax^2+by^2 +f_{\geq 4}(z,t)$  & $2$ & $(0,1,1,1)$&\\
\hline cAx/4 & $ax^2+by^2 +f_{\geq 2}(z,t)$ & $4$ & $(1,3,1,2)$&
$f_2(0,1)=0$\\
\hline cD/2&$x^2+f_{\geq 3}(y,z,t)$ & $2$ &$(1,0,1,1)$&\\
\hline cD/3&$x^2+f_{\geq 3}(y,z,t)$ &$3$ &$(0,2,1,1)$&
$f_3(1,0,0)\neq 0$\\
\hline cE/2&$x^2+y^3+f_{\geq 4}(y,z,t)$ &$2$
&$(1,0,1,1)$&\\
\hline
\end{tabular}
$$
\end{thm} 

\begin{complement} The corresponding quotient singularity 
is terminal iff  the equations   satisfy 2 obvious
conditions:
\begin{enumerate}
\item The equations   define a terminal hypersurface
singularity.

\item The equations   are graded homogeneous.
\end{enumerate} 
 With these assumptions,  a terminal singularity corresponds to exactly
one case on the above list.
\end{complement}

Proof.  By looking at the list of (\ref{hind.overC.thm}), we see
that the assumptions of (\ref{hind.gen.thm}) are satisfied.
Hence we know that $X$ is of the form
$(F^K=0)/\frac1{n}(a_x,a_y,a_z,a_t)$ where 
$\frac1{n}(a_x,a_y,a_z,a_t)$ is on the list of
(\ref{hind.overC.thm}). 

Once we know a $\z_n$-grading on $K[[x,y,z,t]]$ and a
graded homogeneous power series $F^K$, we can try to bring it
to some normal form using the methods (\ref{nf.meth}) and
(\ref{nf.meth.weighted}). They are  set up in such a way
that if $F^K$ is homogeneous in a
$\z_n$-grading the all coordinate changes  respect the
grading.

The proofs of (\ref{ca.thm}, \ref{cd.thm}, \ref{ce.thm}) 
remain unchanged. The only difference is in 
(\ref{morse.lem}). It is not true that a quadratic form can
be diagonalized using a linear transformation which respects
the 
$\z_n$-grading. The best one can achieve is a sum of forms in
disjoint sets of variables $\sum q_i$ where each $q_i$ is
either $au_i^2$ or $u_iv_i$. The latter case is necessary
iff the two variables have different
$\z_n$-grading.

In the $cD$ and $cE$ cases the quadric has rank 1, so it can
be diagonalized.

In the $cA/2$ and $cAx/2$ cases every grade 0 quadric is
diagonalizable.

In the $cAx/4$ case  $x^2,xz,y^2, z^2$ are the only grade 2
quadratic monomials. A quadratic form like this can again be
diagonalized.

Finally let us  look at the $cA/n$-case for $n\geq 3$. The
only grade 0 degree 2 monomials are $xy, t^2$ and  $xz$ if
$r=-1$ or $yz$ if $r=1$. We need to get a rank $\geq 2$
quadric, so $xy $ (or $xz$ if $r=-1$, 
$yz$ if
$r=1$) must appear. In the $r=\pm 1$ case we may need to
perform a linear change of variables to get the normal form
$xy+f(z,t)$. 
\qed

\section{The topology of terminal hypersurface singularities}

Let $0\in X$ be a real singularity.  It's real points
$X(\r)$ form a topological space, which can be triangulated
(cf. \cite[9.2]{BCR87}). We may assume that $0$ is a vertex of the
triangulation. Then locally near
$0$, $X(\r)$ is PL-homeomorphic to the cone over a
simplicial complex 
$L=L(X(\r))$, which is called the {\it link} of
$0$ in $X(\r)$. The  local topology of $X(\r)$ at $0$ is
thus determined by $L$. 

In general one needs to contemplate the dependence of $L$ on
various choices made. I am mainly interested in the case
when $X$ is a 3-dimensional isolated singularity. In this
case $L$ is a compact surface (without boundary) and so $L$
and $X(\r)$ determine each other up to homeomorphism.

The aim of this section is to  classify terminal
singularities over $\r$ according to their local topology.
To be precise, we give a classification in the $cA$ cases
and  provide a procedure in the $cD$ and $cE$ cases which
reduces the 3-dimensional problem to some   questions about
plane curve singularities.

\begin{notation} $M\sim N$ denotes that $M$ and $N$ are
homeomorphic.

 $\uplus$ denotes disjoint union.
$M\uplus rN$ denotes  the disjoint union of $M$ and of $r$ copies
of $N$.

$M_g$ denotes the unique compact, closed and orientable surface of
genus
$g$. 
\end{notation}

  We start with a general lemma.

\begin{lem}\label{orient} Let $X$ be a smooth real
hypersurface. Then
$X(\r)$ is orientable.
\end{lem} 

Proof. Let $X=(f=0)$ be a real equation where $f\in
\r[x_1,\dots,x_n]$ or $f\in \r\{x_1,\dots,x_n\}$. At each
point $p\in X$, $X$ divides a neighborhood of $p$ into two
halves. $f$ is positive on one half and negative on the
other half. Choosing a sign thus determines an orientation.
\qed

\begin{thm}\label{ca1-top.thm}  The following table gives a
complete list of 3-dimensional terminal  singularities of
type $cA_1$ over
$\r$.  

In the table $n\geq 1$.  Case 4, $n=1$ and case 5, $n=1$ are
isomorphic. Aside from this, two  singularities are
isomorphic iff they correspond to the same case and the same
value of $n$.
$$
\begin{tabular}{|c|l|c|}
\hline case & \qquad equation & $L$  \\
\hline $cA_1(1)$ & $x^2+y^2+ z^2\pm t^{2n+1}$ & $S^2$ \\
\hline $cA_1(2)$ & $x^2+y^2- z^2\pm t^{2n+1}$ & $S^2$ \\
\hline $cA_1(3)$& $x^2+y^2+z^2+ t^{2n}$ & $\emptyset$ \\
\hline $cA_1(4)$& $x^2+y^2+z^2-t^{2n}$ & $S^2\uplus S^2$ \\
\hline $cA_1(5)$& $x^2+y^2-z^2+ t^{2n}$ & $S^2\uplus S^2$ \\
\hline $cA_1(6)$& $x^2+y^2-z^2- t^{2n}$ & $S^1\times S^1$ \\
\hline
\end{tabular}
$$
\end{thm}

Proof.  The equations follow from (\ref{ca.thm}), once we
note that
 after multiplying by $\pm 1$  we may assume that  the
quadratic part has at least 2 positive eigenvalues.

The topology is easy to figure out. Since all the claims are
special cases of the next result, I discuss them in more
detail there.\qed

\begin{thm}\label{can-top.thm} A  3-dimensional terminal 
singularity of type $cA_{>1}$ over $\r$ is equivalent to
 a form 
$$ x^2\pm y^2\pm h(z,t)\prod_{i=1}^m f_i(z,t)=0,
$$
 where the $f_i$ are irreducible power series (over $\r$)
such that
$(f_i(z,t)=0)$  changes sign on $\r^2\setminus\{0\}$ and
$h(z,t)$ is  positive on $\r^2\setminus\{0\}$. The
following table gives a complete list of the  possibilites
for the topology of $X(\r)$.
$$
\begin{tabular}{|l|l|c|}
\hline \quad case & \qquad equation & $L(X(\r))$  \\
\hline $cA_{>1}^+(0,+)$& $x^2+y^2+h$ & $\emptyset$ \\
\hline $cA_{>1}^+(0,-)$& $x^2+y^2-h$ & $S^1\times S^1$ \\
\hline $cA_{>1}^+(m)$& $x^2+y^2\pm hf_1\cdots f_m$ &
$\uplus m S^2$\\
\hline $cA_{>1}^-(0)$& $x^2-y^2\pm h$ & $S^2\uplus S^2$ \\
\hline $cA_{>1}^-(m)$& $x^2-y^2\pm hf_1\cdots f_m$  &$M_{m-1}$\\
\hline
\end{tabular}
$$
\end{thm}

Proof. We already have the form $x^2\pm y^2+f(z,t)$ by
(\ref{ca.thm}). Write
$f$ as a product of irreducible power series over $\r$.
Those factors which do not vanish on $\r^2\setminus\{0\}$
are multiplied together to get $h$. By writing
$\pm h$ we may assume that $h$ is   positive on
$\r^2\setminus\{0\}$. (Since the signs of the other factors
are not fixed, the sign of $h$ matters only if there are no
other factors.) Let
$f_i$ be the remaining factors of $f$.

Assume now that we are in the $cA^+$-case: $x^2+ y^2\pm
h\prod f_i$. Projection to the $(z,t)$-plane is a proper map
whose fibers are  as follows:
\begin{enumerate}
\item $S^1$ if $\pm h(z,t)\prod f_i(z,t)<0$, 
\item  a point  if $\pm h(z,t)\prod f_i(z,t)=0$, 
\item empty if $\pm h(z,t)\prod f_i(z,t)>0$.
\end{enumerate}

If $m=0$ then $X(\r)\setminus\{0\}$ is  a circle bundle
over  either
$\r^2\setminus\{0\}$ or over the empty set.  The first case
gives $L\sim S^1\times S^1$ by (\ref{orient}). 

If $m>0$, we have to describe the semi-analytic set
$U:=(\prod f_i(z,t)\leq 0)\subset
\r^2$. Semi-analytic sets can be triangulated
(cf. \cite[9.2]{BCR87}), thus  in a neighborhood of the origin, $U$
is the cone over
$U\cap (z^2+t^2=\epsilon)$. 

Each $(f_i=0)$ is an irreducible curve germ over $\r$, thus
homeomorphic to $\r^1$.  So each $f_i$ has 2 roots on the
circle
$(z^2+t^2=\epsilon)$. Hence $U\cap (z^2+t^2=\epsilon)$ is
the disjoint union of $m$ closed arcs. Therefore $L$ has $m$
connected components, each homeomorphic to $S^2$.

The second possibility is the $cA^-$-case: $x^2- y^2- h\prod
f_i$. (The two choices of $\pm h$  are equivalent by
interchanging
$x$ and $y$.)  Here we project to the $(y,z,t)$-hyperplane.
The fiber  over  a point $(y,z,t)$ is
\begin{enumerate}
\item 2 points if $y^2+ h(z,t)\prod f_i(z,t)>0$, 
\item 1 point if $y^2+ h(z,t)\prod f_i(z,t)=0$, 
\item empty if $y^2+ h(z,t)\prod f_i(z,t)<0$. 
\end{enumerate} Thus we have to determine the region 
$$ U:=\{y^2+ h(z,t)\prod f_i(z,t)\geq 0\}\subset  
(y^2+z^2+t^2=\epsilon)\sim S^2,
$$ and then take its double cover to get $L$. 

If $m=0$ then  $U=S^2$ and so $L=S^2\uplus S^2$. If $m>0$
then $h\prod f_i$ is negative on $m$ disjoint arcs in the
circle
$(z^2+t^2=\epsilon)$, and $y^2+ h(z,t)\prod f_i(z,t)$ is
negative in contractible neighborhoods of these intervals.
Thus
$U=S^2\setminus(\mbox{$m$-discs})$ and so $L$ is a surface
of genus
$m-1$, orientable by (\ref{orient}). 
\qed

\begin{exmp} It is quite instructive to consider the following
incorrect approach to the topology of $cD$ and $cE$-type
singularities. I illustrate it in the $cE_6$-case.

Over $\r$, a surface singularity of type $E_6$ is $x^2+y^2\pm z^4$.
In both cases, projection to the $(x,z)$-plane is a homeomorphism.
Consider    a $cE_6$-type point $X$. If $L(X(\r))$ has several
connected components, then a suitable hyperplane intersects at
least two of them. By a small perturbation we obtain an
$E_6$-singularity as the intersection, thus we conclude that
$L(X(\r))$ is connected. This is especially suggestive if we note
that instead of a plane we could use a small perturbation of any
smooth hypersurface.

Unfortunately the conclusion is false, as we se in
(\ref{ce6.exmp}).   $L(X(\r))$ can have several components,
and some of them are not seen by general hypersurface sections.
These look like very ``thin" cones, as opposed to the main component
which is  ``thick".  It would be interesting to give precise
meaning to this observation and to see its significance in the
study of singularities.
\end{exmp}

The following approach to the   topology of $cD$ and $cE$-type
singularities is taken from
\cite[Sec.12]{AGV85}.

\begin{say}[Deformation to the weighted tangent cone]{\ }

Let $X:=(f(x_1,\dots,x_n)=0)$ be a hypersurface singularity.
For simplicity of notation I assume that $f$ converges for
$|x_i|< 1+\delta$. Assign integral weights to the variables
$w(x_i)=w_i$ and write $f$ as the sum of weighted
homogeneous pieces
$$ f=f_d+f_{d+1}+f_{d+2}+\dots,
$$ where $f_s$ is weighted homogeneous of degree $s$. For a
parameter
$\lambda\neq 0$ set
\begin{eqnarray*} f^{\lambda}(x_1,\dots,x_n):=
\lambda^{-d}f(\lambda^{w_1}x_1,\dots,\lambda^{w_n}x_n)\\ =
f_d+\lambda f_{d+1}+\lambda^2f_{d+2}+\dots
\end{eqnarray*} This suggests that if we define $f^0:=f_d$
then 
$$ X^{\lambda}:=(f^{\lambda}=0)\qtq{for} \lambda\in
\r
$$
 is a ``nice" family of hypersurface singularities. For
$\lambda\neq 0$ they are all isomorphic to $(f=0)$ and for
$\lambda=0$ we obtain the weighted tangent cone
$(f_d=0)$. 

This can be used to determine the topology of $X(\r)$ in 2
steps. First describe $X^0$ and then try to relate
$X^{\lambda}$ and $X^0$ for small values of $\lambda$.

Let $w$ be a common multiple of the $w_i$ and set
$u_i=w/w_i$.
\end{say}

\begin{prop}\label{def.nc.prop}
 Notation as above. Assume that $X=(f=0)$ is an isolated
hypersurface singularity.  

Then there is a $0<\lambda_0$ such that for every 
$0<\lambda\leq \lambda_0$
\begin{enumerate}
\item $L^{\lambda}:=X^{\lambda}\cap (\sum x_i^{2u_i}=1)$ is
smooth, and
\item $X^{\lambda}\cap (\sum x_i^{2u_i}\leq 1)$ is
homeomorphic to the cone over
$L^{\lambda}$.
\end{enumerate}
\end{prop}

Proof. The map $\r^m\to \r^+$ given by
$(x_1,\dots,x_n)\mapsto\sum x_i^{2u_i}$ is proper.  Thus its
restriction to
$X^{\lambda}$ is also proper. The proposition follows once
we establish that the resulting map 
$$ t:X^{\lambda}\to \r
$$ has no critical points with critical value in $(0,1]$ for
$0<\lambda\leq \lambda_0$.

The critical values of a real algebraic morphism form a
semi-algebraic set (cf. \cite[9.5]{BCR87}), thus there is a
$0<\mu_0$ such that
$t:X^1\to \r$ has no critical values in $(0,\mu_0]$.  The
following diagram is commutative
$$
\begin{array}{ccc} X^{\lambda} & \stackrel{x_i\mapsto
\lambda^{-w_i}x_i}{\longrightarrow} & X^1\\ t\downarrow\  &
& \
\downarrow t\\
\r & \stackrel{s\mapsto \lambda^{-w}s}{\longrightarrow}& \r
\end{array}
$$ which shows that (\ref{def.nc.prop}) holds with 
$\lambda_0=\mu_0^{1/w}$.\qed

So far we have not done much, but the advantage of this
approach is that we can view $L^{\lambda}$  as a deformation
of the compact real algebraic variety $L^0$.  If $L^0$ is
smooth then this deformation is locally trivial
differentiably. Thus we obtain:

\begin{cor} Assume that $(f_d=0)$ defines  an isolated
singularity. Then $L^{\lambda}$  is diffeomorphic to
$L^0$.\qed
\end{cor}

This is sufficient to describe the  the topology of 
``general" 
members    of  several  families of terminal
singularities:

\begin{cor}\label{ce.top.gen} Let 
 $X$  be a terminal
singularity given by one of the of  the following   equations:
\begin{enumerate}
\item[$cD_4$] \quad $x^2+f_{\geq 3}(y,z,t)$,  where $f_3=0$
has no real singular point.
\item[$cE_6$] \quad $x^2+y^3+yg_{\geq 3}(z,t)+h_{\geq
4}(z,t)$,  where
$h_4$  has no multiple real linear factor.
\item[$cE_8$] \quad $x^2+y^3+yg_{\geq 4}(z,t)+h_{\geq
5}(z,t)$, where
$h_5$  has no multiple real linear factor.
\end{enumerate} 
Then:
\begin{enumerate}
\item[$cD_4$] \quad  $L(X(\r))\sim S^2$ if $(f_3=0)\subset \r\p^2$
has one connected component and $L(X(\r))\sim S^2\uplus (S^1\times
S^1)$ if
$(f_3=0)$ has two connected components.
\item[$cE_6$]  \quad  $L(X(\r))\sim S^2$.
\item[$cE_8$] \quad  $L(X(\r))\sim S^2$.
\end{enumerate}
\end{cor}

Proof. We use deformation  to the weighted tangent cone
with  (suitable integral multiples of the) weights
$(1/2,1/3,1/3,1/3)$ in the $cD_4$-case, 
$(1/2,1/3,1/4,1/4)$ in the $cE_6$-case, and 
$(1/2,1/3,1/5,1/5)$ in the $cE_8$-case.  The equations for
$X^0$ are
$x^2+f_3(y,z,t)=0$, 
$x^2+y^3+h_4(z,t)=0$  and $x^2+y^3+h_5(z,t)=0$.  Our
conditions guarantee that $X^0(\r)$ has isolated
singularities, thus it is sufficient to determine
$L(X^0(\r))$. 

In the $cE$-cases, projection to the $(x,z,t)$ hyperplane is a
homeomorphism from $X^0(\r)$ to $\r^3$, thus $L(X^0(\r))\sim
S^2$

In the $cD_4$-cases we project to the $(y,z,t)$-hyperplane.
As in the proof of (\ref{can-top.thm}), we can get
$L(X^0(\r))$ once we know the set
$U\subset (y^2+z^2+t^2=1)$ where $f_3$  is nonnegative. The
boundary
$\partial U$ doubly covers the projective curve
$(f_3=0)\subset \r\p^3$. If
 $(f_3=0)\subset \r\p^2$ has one connected component  then
it is a pseudo-line and $\partial U$ is a connected double
cover,  hence $U$ is a disc. If $(f_3=0)$ has two connected
components, then  one is a pseudo-line the other an oval.
$\partial U$ has 3 connected components, and $U$ is a disc
plus an annulus.  Thus $L(X^0(\r))\sim S^2\uplus (S^1\times
S^1)$.
\qed

\begin{rem} In the $cE$ cases of the above example,  projection to
the
$(x,z,t)$ plane is a homeomorphism from $X^0(\r)$ to $\r^3$
even if $h_4$ or $h_5$ have multiple factors. In these
cases, however, we can not conclude that
$X(\r)$ is also homeomorphic to $\r^3$. In fact we   see in
(\ref{ce6.exmp}) that this is not always true.
\end{rem}

Similar arguments work in some of the $cD_{>4}$-cases:

\begin{cor}\label{cd.top.gen} Let $X$  be a terminal
singularity  given by  equation
$x^2+y^2z+h_{\geq s}(z,t)$,   where $z\not\vert h_{s}$ and
$h_{s}$  has no multiple real linear factors. Let $s$ be the
number of real linear factors of $h_{s}$.  There are three
cases:
\begin{enumerate}
\item $s=2r+1$ and 
$L(X(\r))\sim M_r\uplus rS^2$;
\item $s=2r,\ h(0,1)<0$ and 
$L(X(\r))\sim M_r\uplus (r-1)S^2$;
\item $s=2r,\ h(0,1)>0$ and 
$L(X(\r))\sim M_{r-1}\uplus rS^2$.
\end{enumerate}
\end{cor}

Proof. We use deformation  to the weighted tangent cone with
weights
$(1/2,(s-1)/2s,1/s,1/s)$. Thus we need to figure out the
topology of 
$x^2+y^2z+h_{(n-1)}(z,t)=0$. As before, this reduces to
understanding the set where $y^2z+h_{s}(z,t)\leq 0$.  This
can be done by projecting to the $(z,t)$ plane. Details are
left to the reader.
\qed

\begin{say}\label{choose.weights} For many terminal
singularities one can not choose weights so that the
weighted tangent cone has an isolated singularity at the origin, but
in all cases it is possible to choose weights so that the
weighted tangent cone has 1-dimensional singular locus:
$$
\begin{tabular}{|c|l|ccc|}
\hline name &\quad equation is \quad $x^2+$ &$w(y)$ &$w(z)$&$w(t)$\\
\hline $cD_4$& $f_{\geq 3}(y,z,t)$ &
$\frac13$&$\frac13$&$\frac13$\\
\hline $cD_{>4}(1)$&  $y^2z+ayt^r+h_{\geq s}(z,t)$  & 
       $\frac{s-1}{2s}$&$\frac1{s}$&$\frac1{s}$\\
\hline $cD_{>4}(2)$&  $y^2z\pm yt^r+h_{\geq s}(z,t)$  & 
       $\frac{r-1}{2r-1}-$\mbox{\scriptsize $\epsilon$}&
    $\frac1{2r-1}+$\mbox{\scriptsize $2\epsilon$}
&$\frac1{2r-1}+\frac{\epsilon}{r}$\\
\hline $cE_6$ &$y^3+yg_{\geq 3}(z,t)+h_{\geq 4}(z,t)$ &
        $\frac13$&$\frac14$&$\frac14$\\
\hline $cE_7$ &$y^3+yg_{\geq 3}(z,t)+h_{\geq 5}(z,t)$ &
        $\frac13$&$\frac29$&$\frac29$\\
\hline $cE_8$ &$y^3+yg_{\geq 4}(z,t)+h_{\geq 5}(z,t)$ &
        $\frac13$&$\frac15$&$\frac15$\\
\hline
\end{tabular}
$$
In the $cD_{>4}$ case  we asume that $h_s\neq 0$ and use
the first weight sequence if $a=0$ or $2r>s+1$ and  the
second weight sequence if
$2r\leq s+1$, where $\epsilon$ is a  small positive number. 
(We could use $\epsilon=0$ except when $2r=s-1$.) Let
$(w_1,w_2,w_3,w_4)$ be integral multiples of these weights.
The weighted tangent cone, and its singularities are the
following:
$$
\begin{tabular}{|c|l|l|}
\hline name &weighted tangent cone& \qquad singularities  \\
\hline $cD_4$& $x^2+f_3(y,z,t)$ &at singular pts of
$(f_3=0)$\\
\hline $cD_{>4}(1)$&  $x^2+y^2z+h_s(z,t)$&at multiple real
factors of
$zh_s$\\
\hline $cD_{>4}(2)$&  $x^2+y^2z\pm yt^r$ & at the $z$-axis\\
\hline $cE_6$ &$x^2+y^3+h_4(z,t)$& at multiple real factors
of $h_4$\\
\hline $cE_7$ &$x^2+y^3+yg_3(z,t)$& at  real factors of
$g_3$\\
\hline $cE_8$ &$x^2+y^3+h_5(z,t)$& at multiple real factors
of $h_5$\\
\hline
\end{tabular}
$$ These equations have the form $x^2+F(y,z,t)$ and the
deformation to the weighted tangent cone leaves this form
invariant:
$$
x^2+F^{\lambda}(y,z,t)=x^2+\lambda^{-d}F(\lambda^{w_2}y,\lambda^{w_3}z,\lambda^{w_4}t).
$$

Set 
$$ U^{\lambda}:=\{(y,z,t)\vert F^{\lambda}(y,z,t)\leq
0\subset (y^{2u_2}+z^{2u_3}+t^{2u_4}=1)\}.
$$
$U^{\lambda}$ is a semi-algebraic set and its boundary is
the real algebraic curve
$$ C^{\lambda}:=\{(y,z,t)\vert F^{\lambda}(y,z,t)=0\subset
(y^{2u_2}+z^{2u_3}+t^{2u_4}=1)\}.
$$ We have  established the following:
\end{say}

\begin{prop} $C^{\lambda}$ is a deformation of the real
algebraic curve
$C^0$ inside the smooth real algebraic surface 
$(y^{2u_2}+z^{2u_3}+t^{2u_4}=1)$.\qed
\end{prop}

\begin{say} The deformations of real algebraic curves can be
understood in two steps (cf. \cite{Viro90}). Put small discs around
the singularities. Outside the discs all small deformations are
topologically trivial and inside the discs we have a local
problem involving real curve singularities. Here we have the
advantage that
$(y^{2u_2}+z^{2u_3}+t^{2u_4}=1)$ is an affine algebraic
surface, thus we can choose the local deformations
independently and they can always be patched together.

Thus we can describe the possible cases for $C^{\lambda}$,
and thereby the topological types of the corresponding
3-dimensional terminal singularities, if we can describe the
deformations of the occurring real plane curve
singularities. By looking at the equations we see that the
only singularities that we have to deal with are the
2-variable versions of the Du Val singularities:
\begin{enumerate}
\item[$A_n$] \quad $y^2\pm z^{n+1}=0$;
\item[$D_n$] \quad  $y^2z\pm z^{n-1}=0$;
\item[$E_6$]  \quad $y^3+z^4=0$;
\item[$E_7$]  \quad $y^3\pm yz^3=0$;
\item[$E_8$]  \quad $y^3+z^5=0$.
\end{enumerate}

For all of these cases, a complete list of the topological types of
real deformations is known \cite{Chislenko88}.  The list  can also
be found in \cite[Figs. 16--28]{Viro90}, which  contains   many
further examples.
\end{say}

\begin{exmp} Consider for instance the $cD_4$ cases.  The
various possibilities for $f_3$ are easy to enumerate. The
most interesting  is   $f_3=yzt$. Here $C^0$ is the union of
the 3 coordinate hyperplanes intersecting
$(y^{6}+z^{6}+t^{6}=1)$. We have 6 singular points of type
$u^2-v^2=0$. At each of them we can choose a deformation
$u^2-v^2\pm \epsilon=0$. This gives $2^6$ possibilities. The
symetries of the octahedron act on the configurations so it is easy
to get a complete list.

 At the end we get   $7$ possible topological
types for $L(X(\r)$ where
$X=(x^2+yzt+f_{\geq 4}(y,z,t)=0)$:
$$
 M_2, M_1\uplus S^2, M_1, S^2, 2S^2, 3S^2, 4S^2.
$$
 It truns out that these exhaust all the cases given by
$cD_4$.
\end{exmp}

\begin{exmp}\label{ce6.exmp} Consider the $cE_6$-type points
$$
x^2+y^3+yg_{\geq 3}(z,t)\pm z^2t^2+h_{\geq 5}(z,t).
$$
Using the methods of (\ref{nf.meth.weighted}) these can be brought
to the form
$$
x^2+y^3\pm z^2t^2+ya(z)+yb(t)+c(z)+d(t).
$$
The weighted tangent cone, $(x^2+y^3\pm z^2t^2=0)$
 is singular along the $z$ and $t$-axes. In order to understand the
sigularity type of $C^{\lambda}$, say along the positive
$z$-halfaxis, set
$t=\epsilon$. We get an equation
$$
x^2+y^3\pm z^2\epsilon^2+ya(z)+yb(\epsilon)+c(z)+d(\epsilon).
$$
$\mult_0a\geq 3$  and $\mult_0c\geq 5$, thus all the terms
involving $z$ can be absorbed into $z^2$, and we obtain
the equivalent form
$$
x^2\pm z^2+y^3+yb(\epsilon)+d(\epsilon).
$$
The cubic $y^3+yb(\epsilon)+d(\epsilon)$ has 3 real roots if
$4b(\epsilon)^3+27d(\epsilon)^2<0$ and 
1 real root if
$4b(\epsilon)^3+27d(\epsilon)^2>0$.

$C^0$ is homeomorphic to $S^1$ and it has 4 singular points
(along the $z$ and $t$ halfaxes). $C^{\lambda}$ is a smooth curve
which has an oval near a singular point of $C^0$ if the
corresponding cubic has 3 real roots and no ovals if only 1
real root. Thus $C^{\lambda}$ has at most 5 connected components.

We can also determine the location of the ovals relative to the
``main component" of $C^{\lambda}$. In deforming $z^2-y^3=0$,
the oval can appear only toward the negative $y$-direction. 
Putting all this together, we get the following possibilities
for $L(X(\r))$:
$$
\begin{array}{ll}
 rS^2,\ 1\leq r\leq 5 &\mbox{in the $+z^2t^2$-case, and} \\
 M_r,\ 0\leq r\leq 4 &\mbox{in the $-z^2t^2$-case.}
\end{array}
$$
\end{exmp}

\section{The topology of terminal quotient singularities}

Let $0\in X$ be  3-fold terminal singularity and $\pi:\tilde
X\to X$ its index one cover. As we proved,
$X=\tilde X/\frac{1}{n}(a_1,\dots,a_m)$ where $n$ is the
index of $X$ and the $a_i$ are integers. We use this
representation to determine the topology of $X$ in terms of
the already known topology of $\tilde X$.

The main question is to determine the real points of
${\Bbb A}^m/\frac{1}{n}(a_1,\dots,a_m)$.   Let
$\sigma:\c^m\to \c^m$ be the corresponding action of $1\in
\z_n$.

The  answer depends  on the parity of $n$. First we discuss
the odd index cases which are easier.

\begin{prop}\label{odd.quot}
 Assume that $n$ is odd and set
$Y={\Bbb A}^m/\frac{1}{n}(a_1,\dots,a_m)$. Then the induced
map
$\r^n\to Y(\r)$ is a homeomorphism.
\end{prop}

Proof. Let $R\subset \r[[x_1,\dots,x_m]]$ denote the ring of
invariant functions. A point $P\in {\Bbb A}^m$ maps to a
real point of $X$ iff
$f(P)\in \r$ for every $f\in R$. Let $\epsilon$ be a
primitive $n^{th}$ root of unity. If $P=(p_1,\dots,p_m)$ is
real then
$$
\sigma^b(p_1,\dots,p_m)=(\epsilon^{ba_1}p_1,\dots,\epsilon^{ba_m}p_m)
$$
 is also real iff $\sigma^b(P)=P$. This shows that the
quotient map
$\r^n\to X(\r)$ is an injection.

Let $Q\in X(\r)$ be a point. Then $\pi^{-1}(Q)\subset {\Bbb
A}^m$ has an odd number of closed points over $\c$ (usually
$n$ of them) and as a scheme it is defined over $\r$. Thus
it has  a real point, hence
$\r^n\to X(\r)$ is also surjective.\qed

\begin{cor}\label{oddind.termquot} Let $0\in X$ be  a 3-fold
terminal singularity of odd index and $\pi:\tilde X\to X$
its index one cover. Then $\pi: \tilde X(\r)\to X(\r)$ is a
homeomorphism.\qed
\end{cor}

The even index case  is more subtle. For purposes of
induction we allow the case when $n$ is odd. Consider the
action
$\frac{1}{n}(a_1,\dots,a_m)$ on ${\Bbb A}^m$. Write
$n=2^sn'$ where
$n'$ is odd. Let $\eta$ be a  primitive $2^{s+1}$-st root of
unity  and 
$j:\r^m\to \c^m$  the map
$$ j(x_1,\dots,x_m)=(\eta^{a_1}x_1,\dots,\eta^{a_m}x_m).
$$ (If $n$ is odd then $\eta=-1$, hence $j(\r^m)=\r^m$.) 
Write
$a_i=2^ca'_i$ such that $a'_i$ is odd for some $i$. If
$s>c$, let
$\tau:\c^m\to \c^m$ be the $\z_2$-action
$$
\tau(x_1,\dots,x_m)=((-1)^{a'_1}x_1,\dots,(-1)^{a'_m}x_m).
$$ For $s=c$   let $\tau$ be the identity. Note that both
$\r^n$ and
$j(\r^n)$ are $\tau$-invariant.

\begin{prop}\label{2ind.quot} Set
$Y={\Bbb A}^m/\frac{1}{n}(a_1,\dots,a_m)$. Define $j$ and
$\tau$ as above. Then $Y(\r)$ is the quotient of $\r^m\cup
j(\r^m)$ by $\tau$.
\end{prop}

Proof.  The proof is by induction on $m$ and $n$. We can
assume that the action is faithful, that is
$\sum b_ia_i=1$ is solvable in integers. Indeed, for non
faithful actions we get the same quotient from a smaller
group action. The definitions  of  $j$ and $\tau$  are set
up such that they  do not change if we change the group this
way.

Set $F:=\prod_i x_i^{b_i}$. By induction on $m$ we know that
(\ref{2ind.quot}) holds on each coordinate hyperplane. Thus
we have to deal with points $P=(p_1,\dots,p_m)$ such that
each
$p_i\neq 0$. 

Assume that $\pi(P)$ is real. Let $\epsilon$ be a primitive
$(2n)^{th}$ root of unity.
$p_i^n$ is real, hence  
 $p_i=\epsilon^{c_i}\cdot(\mbox{real number})$ for some
$c_i\in \z$.  Thus $F(P)=\epsilon^{c}\cdot(\mbox{real
number})$.
$F(\sigma(P))=\epsilon^2F(P)$, hence by replacing
$P$ by $\sigma^r(P)$ for some $r$  we may assume that
$F(P)\in\r$ or $F(P)\in\eta\cdot\r$.

Assume first that $F(P)$ is real. For each $i$ the function
$F^{n-a_i}x_i$ is invariant, hence has a real value at $P$.
Thus $P\in
\r^n$. If $F(P)\in\eta\cdot\r$ then the same agrument shows
that
$p_i\in \eta^{a_i}\cdot \r$, thus $P\in j(\r^m)$.

This shows that $\r^m\cup j(\r^m)\to Y(\r)$ is surjective.
It is also
$\tau$-invariant.  Finally, if $P=(p_1,\dots,p_m)\in
\r^m\cup j(\r^m)$  then
$$
\sigma^s(P)=(\epsilon^{ba_1}p_1,\dots,\epsilon^{ba_m}p_m)
\in \r^m\cup j(\r^m)
$$ iff 
  $\sigma^s(P)=P$ or $\sigma^s(P)=\tau(P)$. Thus 
$\r^m\cup j(\r^m)\to Y(\r)$ is $2:1$ for $s>c$ and $1:1$ for
$s=c$.
\qed

\begin{defn} 
Let $F\in \r[[x_1,\dots,x_m]$ be a power series, homogeneous
of grade $d$ under the grading $\frac{1}{n}(a_1,\dots,a_m)$.
Let $\eta$ be as above. Define  the {\it companion} $F^c$ of $F$ 
with respect to  the action $\frac{1}{n}(a_1,\dots,a_m)$
by
$$
F^c(x_1,\dots,x_m):=\eta^{-d}F(\eta^{a_1}x_1,\dots,\eta^{a_m}x_m).
$$ Note that $F^c\in \r[[x_1,\dots,x_m]]$. 
\end{defn}

\begin{cor}\label{2ind.termquot} Let $0\in X$ be  a 3-fold
terminal singularity of even index and $\pi:\tilde X\to X$
its index one cover. Then 
$$
L(X(\r)) \sim 
L(\tilde X(\r))/(\tau)\uplus 
L(\tilde X^c(\r))/(\tau).
$$
\end{cor}

Proof.  We use the notation of (\ref{2ind.quot}).  Let $F=0$
be the equation of $\tilde X$ and let
$W:=(F=0)\cap(\r^m\cup j(\r^m))$. Then $X(\r)$ is the
quotient of $W$ by $\tau$. 

$(F=0)\cap\r^m=\tilde X(\r)$. $(F=0)\cap j(\r^m)$ can be
identified with the set of real zeros of 
$F(\eta^{a_1}x_1,\dots,\eta^{a_m}x_m)=0$. (The normalizing
factor
$\eta^{-d}$ does not change the set of zeros.)

In the terminal case the group action is fixed point free on
$(F=0)\setminus\{0\}$, thus $(F=0)\cap\r^m$ and $(F=0)\cap
j(\r^m)$ intersect only at the origin.\qed

As a byproduct we obtain the following:

\begin{cor}\label{h.ind.not.isol} Let $0\in X$ be  a 3-fold
terminal singularity of   index$>1$.  Then
$0\in X(\r)$ is  not an isolated point.
\end{cor}

Proof.  Let $\tilde X$ be the index 1 cover. We are done,
unless
$0\in \tilde X(\r)$ is   an isolated point. This happens
only in cases
$cA/2, cAx/2$ (and maybe for $cAx/4$) where the equation of
$\tilde X$ is
$F=x^2+y^2+f(z,t)$ and $f(z,t)$ is positive on
$\r^2\setminus\{0\}$. 

Let us compute  $F^c$. In the $cA/2$ ase we get
$-x^2-y^2+f(iz,t)$ and this has nontrivial solutions in the
$(x,y,t)$-hyperplane. In the $cAx/2$ case we get
$x^2-y^2+f(iz,it)$ and this has nontrivial solutions in the
$(x,y)$-plane. 

In the $cAx/4$ case already $\tilde X$ has nontrivial
$\r$-points. Indeed, here $f$ has grade 2, thus every
$t$-power in it has an odd exponent. Thus $f(z,t)$ is not 
positive on the $t$-axis.
\qed

\begin{say}[Orientability of index 2 quotients]

We have seen in (\ref{orient}) that every real algebraic hypersurface is
orientable, and so are their quotients by odd order groups  
(\ref{odd.quot}).  With index 2 quotients, the question of
orientablilty is interesting.

Consider a quotient $(F=0)/\frac12(w_1,\dots,w_m)$. 
Let $\sigma$ be the corresponding $\z_2$ action.
We can orient
$X:=(F=0)$ by choosing an orientation of $\r^m$ and at each smooth
point of $X$ we choose the normal vector pointing in the 
direction where $F$ is positive.  $\sigma$ preserves the orientation of
$\r^m$  iff $\sum w_i$ is even. The parity of $w(F)$ determines the sign
in
$\sigma(F)=\pm F$. Thus $\sigma$ preserves the induced orientation of
$X$ iff $w(F)+\sum w_i$ is even. If $w(F)+\sum w_i$ is odd, the induced
orientation is not preserved. If $\sigma$ fixes a connected component
of (the nonsingular part of) $X(\r)$, then the corresponding quotient
is not orientable.  If, however, $\sigma$ only permutes the 
connected components of $X(\r)$ then the quotient is still orientable.
Thus we obtain:

\begin{lem}\label{ind2.orient}
Let $0\in X:=(F=0)/\frac12(w_1,\dots,w_m)$ be an isolated singular
point,
with index one cover $\tilde X$ and companion $\tilde X^c$. Then
$L(X(\r))$ is nonorientable iff $w(F)+\sum w_i$ is odd
and $\sigma$  fixes at least one of the connected  components
of $L(\tilde X(\r))$ or $L(\tilde X^c(\r))$.\qed
\end{lem}
\end{say}

\begin{exmp}[The topology of $cA/2$ points] {\ }

The simplest case is ${\Bbb A}^3/{\textstyle
\frac1{2}(1,1,1)}$. 
 The link of ${\Bbb A}^3$ is the sphere
$(x^2+y^2+z^2=\epsilon^2)$. We act by the antipodal map,
and the quotient is $\r\p^2$. The quotient of the purely
imaginary subspace also gives real points, thus
$L(X(\r))\sim 2\r\p^2$.

In the $cA_{>0}/2$ cases write
  $$
X:=(x^2\pm y^2+f(z,t)=0)/\frac12(1,1,1,0)
$$
with cover $\tilde X:=(x^2\pm y^2+f(z,t)=0)$. 
Only even powers of $z$ occur in  $f(z,t)$, thus we can write
$f(z,t)=G(z^2,t)$. As in (\ref{can-top.thm}) we   factor it as 
$G(z^2,t)=\pm h(z,t)\prod_{i=1}^r f_i(z,t)$.

The companion cover is $\tilde X^c=(x^2\pm y^2-G(-z^2,t)=0)$. We have to
be careful since the product decomposition of $G$ is not preserved. A
factor of $h$ may become indefinite and it can also happen that 
two factors  become conjugate over
$\r$.  Thus we write $-G(-z^2,t)=\pm h'(z,t)\prod_{j=1}^{r'} f'_j(z,t)$.

By (\ref{2ind.termquot}), $L(X(\r))=L(\tilde X(r))/\tau\uplus
L(\tilde X^c(r))/\tau$. In all these cases (\ref{ind2.orient}) shows
that if
$\tau$ fixes a connected component, the quotient is not orientable.
Thus if 
$L(\tilde X(r))$  or $L(\tilde X^c(r))$ is connected, the quotient is
not orientable. This holds in all the $cA^-_{>0}(r>0)$ cases.
In the $cA^-_{>0}(0)$ case the equation is $(x^2=y^2+h)$, and each
halfspace of $(x\neq 0)$ contains a unique connected component.
Thus $\tau$ interchanges the two connected components and the quotient
is orentable.

In the $cA^+_{>0}$ case
$(\pm h(z,t)\prod_{i=1}^r f_i(z,t)=0)$ is negative on 
$r$  connected regions $P_1,\dots,P_r\subset \r^2$, and  
$L(\tilde X(\r))$ consist of $r$ copies of $S^2$, one for
each
$P_j$. The involution $\tau$ fixes the
$t$-axis pointwise, thus if one of the half $t$-axes is
contained in some 
$P_j$, then $\tau$ fixes the $S^2$ over that region. The
other copies of
$S^2$  are interchanged. The same holds for
$\tilde X^c$. 

Along the $t$-axis $G(z^2,t)$ and $-G(-z^2,t)$ have
opposite signs. Thus  among the 4 pairs
$$
\begin{array}{cc} (\mbox{positive half $t$-axis},
G(z^2,t)) & (\mbox{positive half $t$-axis}, -G(-z^2,t))\\
(\mbox{negative half $t$-axis}, G(z^2,t))& (\mbox{negative
half $t$-axis}, -G(-z^2,t))
\end{array}
$$ there are two where the function is negative along the
half axis. 

We obtain the following list of possibilities. (We use the notation
$K_r:=S^2\# r\r\p^2$, thus $K_2$ is the Klein bottle.)

$$
\begin{tabular}{|l|l|c|}
\hline \qquad $\tilde X$ & \qquad $\tilde X^c$ & $L(X(\r))$  \\
\hline $cA_0$& $cA_0$ & $\r\p^2\uplus \r\p^2$ \\
\hline $cA_{>0}^-(r>0)$& $cA_{>0}^-(r'>0)$ & $K_r\uplus K_{r'}$ \\
\hline $cA_{>0}^-(r>0)$& $cA_{>0}^-(0)$ & $K_r\uplus S^2$ \\
\hline $cA_{>0}^-(0)$& $cA_{>0}^-(0)$ & $S^2\uplus S^2$ \\
\hline $cA_{>0}^+(r>0)$& $cA_{>0}^+(r'>0)$ &
       $2\r\p^2\uplus {\textstyle\frac{r+r'-2}{2}}S^2$ \\
\hline $cA_{>0}^+(r>0)$& $cA_{>0}^+(0,+)$ &
       $2\r\p^2\uplus {\textstyle\frac{r-2}{2}}S^2$ \\
\hline $cA_{>0}^+(r>0)$& $cA_{>0}^+(0,-)$ &
       $K_2\uplus {\textstyle\frac{r}{2}}S^2$ \\
\hline $cA_{>0}^+(0,-)$& $cA_{>0}^+(0,+)$ &
       $K_2$ \\
\hline
\end{tabular}
$$
Thus $X(\r)$ is not orientable, except in the fourth case. This
indeed occurs:

 Let $\tilde X:=(x^2-y^2+z^{4m}+t^{2n}=0)$
and 
$X:=\tilde X/{\textstyle
\frac12(1,1,1,0)}$ with companion
$\tilde X^c=(-x^2+y^2+z^{4m}+t^{2n}=0)$. 
 $\tilde X^c\cong \tilde X$ and $L(\tilde X(\r))\sim
S^2\uplus S^2$. $\tau$  interchanges the two
copies of $S^2$. Thus $X(\r)$ is orientable and
$ L(X(\r))\sim S^2\uplus S^2$.
\end{exmp}

\vskip1cm

University of Utah, Salt Lake City UT 84112 

kollar math.utah.edu

\end{document}